\documentclass[secthm]{elsart} 
\usepackage{mathrsfs}
\usepackage{amsmath}
\input{amssymb.sty}
\usepackage{stmaryrd}
 
\DeclareMathOperator{\hks}{HKS}

\begin{document}

\begin{frontmatter}

\title{\bf $p$-Fractals and power series--II.}
\subtitle{Some applications to Hilbert-Kunz  theory}

\author[paul]{Paul Monsky}, 
\ead{monsky@brandeis.edu}
\author[pedro]{Pedro Teixeira}
\address[paul]{Brandeis University, Waltham, MA 02454-9110, USA}
\address[pedro]{Union College, Schenectady, NY 12308-3107, USA}
\journal{Journal of Algebra}
\ead{teixeirp@union.edu}

\begin{abstract}
We use the results of our paper {\em  $p$-Fractals and power series--I\/} (Journal of Algebra  280, 2004, pp. 505--536) to prove the rationality of the Hilbert-Kunz series of a large family of power series, including those of the form $\sum_i f_i(x_i,y_i)$, where the $f_i(x_i,y_i)$ are power series with coefficients in a finite field. The methods are effective, as we illustrate with examples.  In the final section, which can be read independently of the others, we
obtain more precise results for the Hilbert-Kunz  function of the 3 variable
power series $z^D-h(x,y)$.
\end{abstract}

\end{frontmatter}

\section{Introduction}

Throughout this paper $p$ denotes a prime number and $\Bbbk$ a field of characteristic $p$, and the letter  $q$ is  reserved for powers of $p$. For ease of notation we  often denote the list of indeterminates $x_1,\ldots,x_s$ by $\mathbf{x}$,
and  their $q$th powers $x_1^q,\ldots,x_s^q$ by  $\mathbf{x}^q$. 

 The {\bf Hilbert-Kunz function} of $f\in {\Bbbk\llbracket \mathbf{x} \rrbracket}$ is the function 
\begin{equation}
\begin{array}{rcl}
n\longmapsto e_n(f)&=&\deg(\mathbf{x}^{p^n},f)\\
&=&\dim_\Bbbk {\Bbbk\llbracket \mathbf{x} \rrbracket}/(\mathbf{x}^{p^n},f).
\end{array}\nonumber
\end{equation}
This is the Hilbert-Kunz function of the  local ring  ${\Bbbk\llbracket \mathbf{x} \rrbracket}/(f)$ with respect to its maximal ideal.      The {\bf Hilbert-Kunz series}  of $f$ is the associated power series $\hks(f)=\sum_{n=0}^\infty e_n(f){\mathsf{z}}^n $, and  the {\bf Hilbert-Kunz multiplicity}  of $f$ is the limit of $e_n(f)/p^{(s-1)n}$ as $n\to \infty$. 
 
  In many cases (see for example  \cite{ch2}, \cite{paul2}, \cite{paul3} and  \cite{tese})  the Hilbert-Kunz series is  a quotient of polynomials with integer coefficients.  (Computer calculations strongly suggest that this is false in
general---certain simple 5 variable $f$ are likely counterexamples.)  In this article we use the ideas and results of \cite{pfractals1} to extend these rationality results to  a large family of power series, including those of the form $\sum_{i=1}^n f_i(x_i,y_i)$, where the $f_i(x_i,y_i)$ are power series with coefficients in a finite field.   (When $n=1$ this is proved
rather simply in \cite{paul1} without the finiteness restriction on $\Bbbk$.)

Let $\mathscr{I}$ be the set of rational numbers in $[0,1]$ whose denominators are powers of $p$. We are  concerned with the following family of functions $\mathscr{I}\to\Qset$:

\begin{defn}\label{def_phi_f}{\em 
Suppose  $f\in \Bbbk\llbracket x_1,\ldots,x_s\rrbracket$. Then $\varphi_f:\mathscr{I} \to \Qset$ is the function $a/q \mapsto q^{-s}\deg(\mathbf{x}^q,f^a)$. 
}\end{defn}

In the notation introduced in  \cite{pfractals1}, 
$\varphi_f$ is the function $\varphi_I$ with $I=(x_1,\ldots,x_s)$, $r=1$ and  $h_1=f$. In that paper we introduced the concept of {\bf $p$-fractals}:  a function $\varphi:\mathscr{I}\to\Qset$ is a  $p$-fractal if the $\Qset$-subspace of $\Qset^\mathscr{I}$ spanned by the functions $T_{q|b}\varphi: t\mapsto \varphi((t+b)/q)$ ($q$ a power of $p$; $b\in \Zset$ with $0\le b<q$) is finite dimensional. Equivalently, $\varphi$ is a $p$-fractal if it is contained in a finite dimensional $\Qset$-subspace of $\Qset^\mathscr{I}$ which is stable under the operators $T_{p|b}$ ($0\le b<p$). We call such a subspace  {\bf $p$-stable}. 
 Theorem 1 of \cite{pfractals1} shows  that all the functions $\varphi_I$ (so in particular the $\varphi_f$) are $p$-fractals when $s=2$ and $\Bbbk$ is finite.

The  following simple proposition establishes the connection between $p$-fractals and the rationality of the Hilbert-Kunz series: 

\begin{prop}\label{rationality}
Let $f\in \Bbbk\llbracket x_1,\ldots,x_s \rrbracket$, and suppose the function $\varphi_f$ is a $p$-fractal.  Then the Hilbert-Kunz  series of $f$ is rational.  
\end{prop}

\begin{pf}
Let $\varphi=\varphi_f$ and $S=T_{p|0}$. Since  $\varphi$ is  a $p$-fractal, the $\Qset$-subspace $V$ of $\Qset^\mathscr{I}$ spanned by  $S^n(\varphi)$, $n=0,1,2,\ldots$,  is finite dimensional, and $S$ maps $V\to V$. So 
the restriction of $S$ to $V$ satisfies  a polynomial  identity of the form
$S^{l}=c_{1}S^{l-1}+c_{2}S^{l-2}+\cdots+c_{l}I$, where $I$ is the 
identity map and $c_i\in \Qset$.   This equation  can be applied to any $S^n(\varphi)$, and shows that the sequence
 $(S^n(\varphi))$ is linearly recurrent. Evaluating at 1, we get a linear recursion 
for $S^n(\varphi)(1)=\varphi(1/p^n)=p^{-sn}e_n(f)$, and the proposition  follows easily. 
\qed
\end{pf}

In the light of Proposition \ref{rationality}, we  define: 

\begin{defn}\label{defn-strongly-rational}{\em 
 A power series $f\in \Bbbk\llbracket \mathbf{x}\rrbracket$  is  {\bf strongly rational} if  $\varphi_f$ is a $p$-fractal.
}\end{defn}

We shall  prove   that strong rationality is preserved under certain   operations: if $f \in \Bbbk\llbracket x_1,\ldots,x_s\rrbracket$ and $g\in \Bbbk\llbracket y_1,\ldots, y_{s'} \rrbracket$ 
 are strongly rational, then so are 
 the product $fg$, 
 the sum $f+g$ and  powers of these power series.   Using    \cite[Theorem 1]{pfractals1} as a starting point, these results  show, for example,  the rationality of the Hilbert-Kunz series of  power series of the form $\sum_i f_i(x_i,y_i)$ and $\sum_i\prod_j f_{ij}(x_{ij},y_{ij})$, where  $f_i(x_i,y_i)$ and $f_{ij}(x_{ij},y_{ij})$ are power series with coefficients in a finite field. 
 
While the proofs that strong rationality is preserved by powers and  products, in the sense described above,  are  relatively straightforward, the proof of the analogous result for sums will require some work---in particular, it will be necessary  to use  the {\em representation ring\/}  introduced by Han and the first author in \cite{paul2}.  
Section \ref{representation_ring} summarizes a few definitions and results of \cite{paul2}, and introduces an endomorphism $\theta$ of the representation ring. In Section \ref{sequences} we attach to each function $\varphi:\mathscr{I}\to \Qset$ a certain sequence of elements of the representation ring, and find conditions on that sequence which are necessary and sufficient for $\varphi$ to be a $p$-fractal. In Section \ref{operations} we use these
conditions to show that if $f(\mathbf{x})$ and $g(\mathbf{y})$ are strongly rational, then the same
is true for $f(\mathbf{x})+g(\mathbf{y})$. We prove analogous results for products and powers;
as we have noted, this together with  \cite[Theorem 1]{pfractals1} gives our rationality
results.

In Section \ref{examples} we show how our techniques can be used to effectively calculate the Hilbert-Kunz series of various power series. We conclude with an  analysis of the Hilbert-Kunz function of $z^D-h(x,y)$, in Section \ref{application}.

\section{The representation ring and the endomorphism $\theta$}\label{representation_ring}

We  summarize  here some of the definitions and results of \cite{paul2}, for the reader's convenience. 
A {\bf $\Bbbk$-object} is a finitely 
generated $\Bbbk[T]$-module on which $T$ acts nilpotently. 
$\Gamma$ is the 
Grothendieck group of
the semigroup of isomorphism classes of $\Bbbk$-objects under 
the usual direct sum. We  introduce a product on $\Gamma$ as follows: 
if two elements of $\Gamma$ are represented by the $\Bbbk$-objects
$M$ and $N$, then their product is the image in $\Gamma$ of the $\Bbbk$-object
$M\otimes_{\Bbbk}N$, where $T$ acts distributively; namely $T(m\otimes 
n)=(Tm)\otimes n +m\otimes (Tn)$.
Now $\Gamma$ endowed with this product  is a 
commutative ring, called the {\bf representation ring}. The zero 
and unity  of $\Gamma$ are respectively the images of the 
zero module and $\Bbbk[T]/(T)$ in $\Gamma$.

For any nonnegative integer $n$,  $\delta_{n}$ is the image of 
$M_{n}=\Bbbk[T]/(T^{n})$ 
in $\Gamma$ (so in 
particular $\delta_{0}=0$ and $\delta_{1}=1$). 
The theory of modules over principal ideal domains  shows that $(\Gamma,+)$ is a 
free abelian group with basis $\{\delta_{1},\delta_{2},\ldots\}$.  In what follows we shall mostly use  a  second basis $\{\lambda_0,\lambda_1,\ldots\}$, where  $\lambda_{n}=(-1)^{n}(\delta_{n+1}-\delta_{n})$. 

The $\lambda_{i}$-coordinate of the image in  $\Gamma$  
of  a $\Bbbk$-object $M$ is
$(-1)^{i}\dim_{\Bbbk} 
T^{i}M/T^{i+1}M$, or $(-1)^{i}(\dim_{\Bbbk} M/T^{i+1}M-\dim_{\Bbbk} 
M/T^{i}M).$ 
This is immediately seen for $M_n=\Bbbk[T]/(T^{n})$, whose image in $\Gamma$ is $\delta_{n}=\sum_{i<n}(-1)^{i}\lambda_{i}$, and evidently  the result  extends additively to arbitrary $\Bbbk$-objects.

\begin{defn}\label{defn-[f]n}{\em 
Let $f\in {\Bbbk\llbracket \mathbf{x} \rrbracket}$. 
Then $\langle f\rangle_{n}$ is the image in $\Gamma$ 
of the $\Bbbk$-object ${\Bbbk\llbracket \mathbf{x} \rrbracket}/(\mathbf{x}^{p^{n}})$,
where $T$ operates by multiplication 
by $f$.  
}\end{defn}

In terms of the $\lambda_{i}$,
\begin{equation}\label{coord-of-[f]n}
\langle f\rangle_{n}
=\sum_{i=0}^{q-1}(\deg(\mathbf{x}^{q},f^{i+1})-
\deg(\mathbf{x}^{q},f^{i}))(-1)^{i}\lambda_{i},
\end{equation}
where $q=p^n$.

If $g$ is another power series in a different set of variables, say 
$g\in \Bbbk\llbracket \mathbf{y} \rrbracket$, then the product $\langle f\rangle_{n}\langle g\rangle_{n}$ is represented by the $\Bbbk$-object $\Bbbk\llbracket \mathbf{x},\mathbf{y}\rrbracket/(\mathbf{x}^q,\mathbf{y}^q)$, where $T$ operates by multiplication 
by $f+g$. So
$\langle f\rangle_{n}\langle g\rangle_{n}=\langle f+g\rangle_{n}$. 

\begin{defn}\label{defn-alpha}{\em 
$\alpha:\Gamma\rightarrow\Zset$ is the $\Zset$-linear map $\sum_{i=0}^n c_i\lambda_i\mapsto c_0$. 
}\end{defn}

Note that   $\alpha(\langle f\rangle_{n})=\deg(\mathbf{x}^{p^{n}},f)=e_{n}(f)$, for all $f\in {\Bbbk\llbracket \mathbf{x} \rrbracket}$. Theorem 1.10 of \cite{paul2} shows that  $\alpha(\lambda_i\lambda_j)=\delta_{i,j}$ ($=1$ if $i=j$; 0 otherwise).

We shall  use  the following multiplication formulas, from Lemma 3.3 and Theorems 2.5, 3.4 and 3.10 of 
\cite{paul2}:
\begin{equation}\label{eq5}
\delta_{i}\delta_{qj}  = i\cdot \delta_{qj}  \ \ \ (0\le i\le q, \ j\ge 1); 
\end{equation}
\begin{equation}\label{eq1}
\lambda_{i}\lambda_{qj}  = \lambda_{qj+i}  \ \ \ (0\le i< q); 
\end{equation}
\begin{equation}\label{eq2}
 \lambda_{i}\lambda_{qj-1} =  
\lambda_{qj-1-i} \ \ \ (0\le i< q);
\end{equation}
\begin{equation}\label{eq3}
\lambda_{q}\lambda_{qj} =  
\lambda_{q(j-1)}+\lambda_{q(j+1)-1}+\lambda_{q(j+1)}, \mbox{\ if $p\not \hskip -.012in |\  j$ and $p\not\hskip -.012in   |\  j+1$;}
\end{equation}
\begin{equation}\label{eq4}
\lambda_{i}\lambda_{j}=\sum_{k=j-i}^{\min(i+j,2p-2-i-j)}\lambda_{k}, \mbox{\ if $i\le j <p$}.
\end{equation}

In the remainder of this  section we introduce a linear map $\theta:\Gamma\to \Gamma$, and show that $\theta$ is actually a ring homomorphism---a  result essential to  the next section. The existence of this endomorphism  is an important
result in its own right, which helps in understanding the structure of the  representation ring.

\begin{defn}\label{defn-phi}{\em 
$\theta$ is the linear operator on $\Gamma$ whose values on the basis $\{\lambda_{i}\}$ are defined as 
follows:
\begin{displaymath}
\theta(\lambda_i)=
\begin{cases}
\lambda_{pi}  & \mbox{if $i$ is even;}\\
\lambda_{pi+p-1} &\mbox{if $i$ is odd.} \\
\end{cases}
\end{displaymath}
}\end{defn}

In what follows, let $n$ be a fixed positive integer and $q=p^n$. 

\begin{lem}
$\theta^{n}(\Gamma)$ is a subring of $\Gamma$.
\end{lem}

\begin{pf}
$\theta^{n}(\Gamma)$ is the $E_{q}$ in Definition 3.12 of  \cite{paul2}. So the lemma follows immediately from Corollary 3.16 of that same paper.  
\qed
\end{pf}

\begin{defn}\label{defn-gamma-n}{\em 
$\Gamma_{n}$ is the additive subgroup of $\Gamma$ generated by the 
$\lambda_{i}$ with $ i<q=p^n$. 
}\end{defn}

\begin{rem}\label{gamma_n_is_subring}{\em
Theorem 3.2 of \cite{paul2} shows that   $\Gamma_{n}$ is   a subring of $\Gamma$.
}\end{rem}

\begin{lem}\label{decompose-gamma}
$\theta^{n}(\Gamma)\Gamma_{n}=\Gamma$.
\end{lem}

\begin{pf} It suffices to show that 
$\lambda_{i}\in\theta^{n}(\Gamma)\Gamma_{n}$, for all $i$. Write 
$i=aq+b$, with $0\le b <q$. 
If $a$ is even, Eq. (\ref{eq1}) shows that 
$\lambda_{i}=\lambda_{aq}\lambda_{b}=\theta^{n}(\lambda_{a})\lambda_{b}\in
\theta^{n}(\Gamma)\Gamma_{n}$, while if $a$ is odd,  Eq. (\ref{eq2}) gives 
 $\lambda_{i}=\lambda_{q(a+1)-1}\lambda_{q-1-b}=
\theta^{n}(\lambda_{a})\lambda_{q-1-b}\in \theta^{n}(\Gamma)\Gamma_{n}.$
\qed
\end{pf}

\begin{lem}\label{case1}
If $u\in \theta^{n}(\Gamma)$ and $v\in \Gamma_{n}$, then 
$\theta(uv)=\theta(u)\theta(v)$.
\end{lem}

\begin{pf}
It is enough to verify that    $\theta(\lambda_{2qi}\lambda_{j})=
\theta(\lambda_{2qi})\theta(\lambda_{j})$ and 
$\theta(\lambda_{2qi-1}\lambda_{j})=
\theta(\lambda_{2qi-1})\theta(\lambda_{j})$, for any $i$  and $j$ with   
$j<q$. 
Those follow from   simple calculations using   equations  (\ref{eq1}) and  (\ref{eq2}). 
\qed\end{pf}

\begin{lem}\label{leminha2}
Suppose $1\le j\le p-2$. Then
\begin{itemize}
\item
$\lambda_{2q-1}\lambda_{q(j+1)-1}=\lambda_{q(j-1)}+\lambda_{q(j+1)-1}+\lambda_{q(j+1)};$
\item
$\lambda_{2q-1}\lambda_{qj}=\lambda_{qj-1}+\lambda_{qj}+\lambda_{q(j+2)-1}$.
\end{itemize}
\end{lem}

\begin{pf}
Equations  (\ref{eq1}) and   (\ref{eq2})   show that $\lambda_{2q-1}\lambda_{q(j+1)-1}=(\lambda_{q-1}\lambda_{q})(\lambda_{q-1}\lambda_{qj})=\lambda_{q}\lambda_{qj}$, so  the first identity  follows from Eq. (\ref{eq3}). 
Also
$\lambda_{2q-1}\lambda_{qj}=\lambda_{q-1}\lambda_{q}\lambda_{qj}=
\lambda_{q-1}(\lambda_{q(j-1)}+\lambda_{q-1}\lambda_{qj}+\lambda_{q(j+1)})$, and  
equations  (\ref{eq1}) and   (\ref{eq2})  give the second identity.
\qed\end{pf}

In what follows,  let  $\mu_i=\theta^n(\lambda_i)$.

\begin{lem}\label{leminha3}
\begin{itemize}
\item If $1\le j \le p-2$, 
$\mu_{1}\mu_{j}=\mu_{j-1}+\mu_{j}+\mu_{j+1}$.  
\item $\mu_{1}\mu_{p-1}=\mu_{p-2}$. 
\end{itemize}
\end{lem}

\begin{pf}
The first identity  is simply a reformulation of  Lemma \ref{leminha2}. 
When $p=2$, 
the second identity  becomes  $(\lambda_{2q-1})^{2}=\lambda_{0}$, which is a special case of Eq. (\ref{eq2}). On the other 
hand, when $p\ne 2$, 
$\lambda_{2q-1}\lambda_{q(p-1)}=\lambda_{q}\lambda_{q-1}\lambda_{q(p-1)}=\lambda_{q}\lambda_{qp-1}=\lambda_{q(p-1)-1}$; this is precisely the second identity.
\qed\end{pf}

\begin{lem}If $i\le j <p$, then
\begin{equation}
\mu_{i}\mu_{j}=\sum_{k=j-i}^{\min(i+j,2p-2-i-j)}\mu_{k}.\nonumber
\end{equation}
\end{lem}

\begin{pf}
It suffices to prove the identity when $i+j<p$ and to show that 
$\mu_{i}\mu_{p-1}=\mu_{p-1-i}$, for $i<p$. Both of these results are proved 
by induction on $i$---Lemma \ref{leminha3} gives the result for $i=1$, and the inductive step relies on  the identity $\mu_i=\mu_1 \mu_{i-1}-\mu_{i-1}-\mu_{i-2}$.  The 
calculations are identical to  those of Lemma 2.4 and Theorem 2.5 of \cite{paul2}.
\qed\end{pf}

Comparing the above lemma and Eq. (\ref{eq4}) we conclude the following:

\begin{lem}\label{case2}
Suppose $i,j <p$. Then 
$\theta^{n}(\lambda_{i}\lambda_{j})=\theta^{n}(\lambda_{i})\theta^{n}(\lambda_{j})$.   
\end{lem}

\begin{thm}\label{phi-is-homo}
$\theta$ is a ring homomorphism.
\end{thm}

\begin{pf} We will show by induction on $n$ that 
$\theta(uv)=\theta(u)\theta(v)$ for all $u$ and $v$ in $\Gamma_{n}$. 
Lemma \ref{case2} gives the result for $n=1$, so  
now suppose  the assertion  holds for 
$n\ge 1$, and let $u,v\in \Gamma_{n+1}$. We may 
assume   that $u=\lambda_{a}$ and 
$v=\lambda_{b}$, with $ a,b <p^{n+1}$. Then as in the proof of Lemma
\ref{decompose-gamma} we 
may write $\lambda_{a}=\theta^{n}(\lambda_{i})\lambda_{j}$ with 
$  i <p$ and $  j< p^{n}$. Similarly, we write 
$\lambda_{b}=\theta^{n}(\lambda_{k})\lambda_{l}$ with 
$  k <p$ and $  l < p^{n}$. Then
$\theta(uv)  =  
\theta\big(\theta^{n}(\lambda_{i})\theta^{n}(\lambda_{k})
\lambda_{j}\lambda_{l}\big)$. 
Since both $\theta^{n}(\Gamma)$ and $\Gamma_{n}$ are closed 
under multiplication, Lemma \ref{case1} shows that this is equal to 
$\theta\big(\theta^{n}(\lambda_{i})\theta^{n}(\lambda_{k})\big)\cdot
\theta(\lambda_{j}\lambda_{l})$. From Lemma \ref{case2} and the 
induction assumption  it follows that $\theta(uv)=
\theta(\theta^{n}(\lambda_{i}))\cdot
\theta(\lambda_{j})\cdot\theta(\theta^{n}(\lambda_{k}))\cdot
\theta(\lambda_{l})$, and one more application of Lemma \ref{case1} 
concludes the proof.
\qed\end{pf}

\section{Coherent sequences}\label{sequences}

In this section we associate to each function $\varphi:\mathscr{I}\to \Qset$ a sequence of elements of the representation ring, and find  necessary and sufficient conditions on that sequence in order for   $\varphi$ to be a $p$-fractal. For our purposes  we will  need to allow rational coefficients, working  with $\Gamma_\Qset=\Gamma\otimes_\Zset \Qset$ rather than $\Gamma$.  The linear map $\alpha\otimes_\Zset 1:\Gamma_\Qset\to \Qset$ will  also be denoted by $\alpha$, by abuse of notation. 

\begin{defn}{\em 
$\Lambda$ is the ring  $(\Gamma_\Qset)^\Nset$ of sequences $\mathbf{u}=(u_{0},u_{1},u_{2},\ldots)$ with entries in $\Gamma_\Qset$. 
}\end{defn}

\begin{defn}\label{defn-L}{\em 
Given any function $\varphi:\mathscr{I}\to \Qset$, $\mathscr{L}(\varphi)$ is the element of $\Lambda$ whose $n$th entry is
\begin{equation}
\mathscr{L}_n(\varphi)=\sum_{i=0}^{q-1}\textstyle\left(\varphi\left(\frac{i+1}{q}\right)-\varphi\left(\frac{i}{q}\right)\right)(-1)^i\lambda_i,\nonumber
\end{equation}
where $q=p^n$.  $\mathscr{L}:\Qset^\mathscr{I}\to \Lambda$ is the map $\varphi\mapsto\mathscr{L}(\varphi)$.
}\end{defn}

In particular, if $f\in \Bbbk\llbracket x_1,\ldots,x_s \rrbracket$ and  $\varphi=\varphi_f$ (see Definition \ref{def_phi_f}) then comparing the above definition to Eq. (\ref{coord-of-[f]n}) we see that   $\mathscr{L}_n(\varphi_f)=p^{-sn}\langle f\rangle_n$. Consequently,  $\alpha(\mathscr{L}_n(\varphi_f))=p^{-sn}e_n(f)$. As noted after Eq. (\ref{coord-of-[f]n}), if $g$ is a power series in the indeterminates $y_1,\ldots,y_{s'} $, then $\langle f+g\rangle_n=\langle f\rangle_n\langle g\rangle_n$, and it follows at once that  $\mathscr{L}(\varphi_{f+g})=\mathscr{L}(\varphi_{f})\mathscr{L}(\varphi_{g}).$

Our first step towards finding the  alternative characterization of $p$-fractals will be to describe the image of the map $\mathscr{L}$.  
Theorem 4.1 of \cite{paul2} introduces an endomorphism  $\psi_p$ of $\Gamma$ with the following property: 
\begin{displaymath}
\psi_p(\delta_{pr+k})=(p-k)\delta_r+k\delta_{r+1},
\end{displaymath}
for all $r$ and $k$ with $0\le k\le p$.
It follows that 
\begin{equation}\label{prop-psi}
\psi_p(\lambda_{pr+k})=(-1)^{pr+k+r}\lambda_r,
\end{equation}
for all $r$ and $k$ with $0\le k< p$. On the $\Bbbk[T]$-module level, $\psi_p$ of a $\Bbbk$-object has the same underlying module, but the new action of $T$ is the $p$th power of the old action.

\begin{defn}{\em 
$\psi$ is the endomorphism $\psi_p\otimes_\Zset 1$ of $\Gamma_\Qset$. 
}\end{defn}

One can easily  see that 
\begin{equation}\label{map-well-defn}
\alpha(\delta_{pi}u)=\alpha(\delta_i\psi(u)),
\end{equation} for all $i\in \Nset$ and $u\in \Gamma_\Qset$. It suffices to verify the formula  for $u=\lambda_{pr+k}$, with $0\le k< p$. Then,   since   $\delta_i=\sum_{j<i}(-1)^j\lambda_j$ and $\alpha(\lambda_i\lambda_j)=\delta_{i,j}$,  both sides of  (\ref{map-well-defn})  are  zero if $i\le r$, and $(-1)^{pr+k}$ otherwise.

\begin{defn}\label{defn-coherent}{\em 
A sequence $\mathbf{u}=(u_0,u_1,u_2,\ldots)$ of $\Lambda$  is {\bf coherent}  if the following properties hold:
 \begin{enumerate}
\item  Each $u_n$ is a linear combination of $\lambda_i$ with $i<p^n$ (i.e., $u_n\in \Gamma_n\otimes_\Zset \Qset$);
\item $\psi(u_{n+1})=u_n$, for all $n$. 
\end{enumerate}
$\Lambda_0$ is the $\Qset$-subspace of $\Lambda$ consisting of the coherent  sequences. 
}\end{defn}

Remark \ref{gamma_n_is_subring} and the fact   that  $\psi$ is a ring homomorphism show that  $\Lambda_0$ is actually  a $\Qset$-subalgebra of $\Lambda$.   A   calculation using  Eq. (\ref{prop-psi}) shows that  $\mathscr{L}(\varphi)$ is coherent, for any $\varphi\in \Qset^\mathscr{I}$.  (For $\varphi=\varphi_f$ this is more conceptually seen using the description of $\psi_p$ on the $\Bbbk[T]$-module level.)
In fact,  $\Lambda_0$  consists {\em precisely\/} of the   $\mathscr{L}(\varphi)$, as  the following lemma shows. 

\begin{lem}\label{onto-coherent}
 $\mathscr{L}:\Qset^\mathscr{I}\to \Lambda$ maps $\Qset^\mathscr{I}$ onto $\Lambda_0$, and\/ $\ker \mathscr{L}$ is the 1-dimensional  subspace of constant functions.
\end{lem}

\begin{pf} Suppose $\mathbf{u}=(u_0,u_1,u_2,\ldots)$ is a coherent sequence. Let  $\varphi:\mathscr{I}\to \Qset$ be the map $i/p^n\mapsto \alpha(\delta_i u_n)$. Note that Eq.  (\ref{map-well-defn}) and the fact that $\mathbf{u}$ is coherent imply that $\varphi$ is well-defined. If we write $u_n=\sum_j (-1)^j a_j \lambda_j$, then $\alpha(\delta_i u_n)=a_0+\cdots+a_{i-1}$.  So $\varphi((i+1)/p^n)-\varphi(i/p^n)=\alpha(\delta_{i+1} u_n)-\alpha(\delta_i u_n)=a_i$, and it follows that $\mathscr{L}(\varphi)=\mathbf{u}$, showing that $\mathscr{L}$ maps $\Qset^\mathscr{I}$ onto $\Lambda_0$. 

To conclude the proof, note that $\mathscr{L}$ maps a function $\varphi$ to 0 if and only if $\varphi((i+1)/p^n)-\varphi(i/p^n)=0$ for all $n$ and $i<p^n$, so the kernel of $\mathscr{L}$ consists of the constant functions.  
\qed
\end{pf}

We give  $\Lambda$ a  structure of $\Gamma$-module, introducing  a product $\Gamma\times \Lambda\to \Lambda$ as follows:
\begin{defn} \label{defn-product} {\em 
Suppose $w\in \Gamma$ and $\mathbf{u}=(u_{0},u_{1},u_{2},\ldots)\in \Lambda$. Then 
\begin{displaymath}
w\cdot\mathbf{u}=(wu_{0},\theta(w)u_{1},\theta^2(w) u_{2},\ldots).
\end{displaymath}
}\end{defn}

\begin{defn}{\em 
$R:\Lambda\to \Lambda$ is the $\Gamma$-linear map    taking  $(u_0,u_1,u_2,\ldots)$ to $(v_0,v_1,v_2,\ldots)$, where $v_n=(-1)^{p^n} \lambda_{p^n-1} u_n.$
}\end{defn}

One sees directly from the above definition and Eq. (\ref{eq2}) that if $\varphi$ is a function $\mathscr{I}\to \Qset$ and $\bar{\varphi}$ is its ``reflection'' $t\mapsto \varphi(1-t)$, then 
\begin{equation}\label{R_and_reflections}
\mathscr{L}(\bar{\varphi})=R(\mathscr{L}(\varphi)).
\end{equation} 
In particular, Lemma \ref{onto-coherent} then  shows that $R$ stabilizes $\Lambda_0$.

\begin{defn}\label{defn-shift}{\em 
$S:\Lambda\to \Lambda$ is the map $(u_{0},u_{1},u_{2},\ldots)\mapsto (u_{1},u_{2},u_{3},\ldots)$. 
}\end{defn}

We shall  refer to $S$ as the {\bf shift operator}.  Quick calculations using Eq. (\ref{eq1}) show the following identities relating $R$ and $S$:
\begin{equation}\label{R_and_S}
S(R(\mathbf{u}))=\lambda_1 S(\mathbf{u}) \mbox{ if $p=2$;}
 \hskip .3in  S(R(\mathbf{u})) =\lambda_{p-1} R(S(\mathbf{u}))  \mbox{ if $p>2$.}
\end{equation}

We now describe the action of the shift operator  on $\mathscr{L}(\varphi)$:

\begin{lem}\label{shift_op_on_Lphi}
Suppose $\varphi$ is a function $\mathscr{I}\to \Qset$. Then  the shift operator acts on $\mathscr{L}(\varphi)$ as follows:
\begin{displaymath}S(\mathscr{L}(\varphi))=\sum_{{k\ \hbox{\tiny  even}}\atop {0\le k<p}}\lambda_k\mathscr{L}(T_{p|k}\varphi)+
\sum_{{k\ \hbox{\tiny  odd}}\atop {0\le k<p}}\lambda_k\mathscr{L}\left(\overline{T_{p|k}\varphi}\right),   
\end{displaymath}
where $\overline{T_{p|k}\varphi}$ denotes the ``reflection''  of ${T_{p|k}\varphi}$, namely the map  $t\mapsto T_{p|k}\varphi(1-t)$. 
\end{lem}

\begin{pf}
Fix a nonnegative integer $n$ and  let $q=p^n$. Then  
\begin{displaymath}
\mathscr{L}_{n+1}(\varphi)=\sum_{i=0}^{pq-1}\textstyle\left(\varphi\left(\frac{i+1}{pq}\right)-
\varphi\left(\frac{i}{pq}\right)\right)(-1)^{i}\lambda_{i},
\end{displaymath} 
and this sum can be split 
up into $p$ sums
\begin{displaymath}
s_{k}=\sum_{i=kq}^{kq+q-1}\textstyle\left(\varphi
\left(\frac{i+1}{pq}\right)-
\varphi\left(\frac{i}{pq}\right)\right)(-1)^{i}\lambda_{i}\ \hskip .5in (0\le k<p).
\end{displaymath}
Suppose first that $k$ is even. Shifting indices and  using  Eq. (\ref{eq1}) we obtain
\begin{eqnarray}
s_{k}&=&  \sum_{i=0}^{q-1}\textstyle\left(\varphi\left(\frac{i+kq+1}{pq}\right)-
\varphi\left(\frac{i+kq}{pq}\right)\right)(-1)^{i}\lambda_{i+kq}\nonumber \\
&=&
\lambda_{kq}\sum_{i=0}^{q-1}\textstyle\left(\varphi\left(\frac{i+kq+1}{pq}\right)-
\varphi\left(\frac{i+kq}{pq}\right)\right)(-1)^{i}\lambda_{i},\nonumber
\end{eqnarray}
which can be rewritten as 
\begin{eqnarray}
s_{k}&=& \lambda_{kq}
\sum_{i=0}^{q-1}\textstyle\left(T_{p|k}\varphi\left(\frac{i+1}{q}\right)-
T_{p|k}\varphi\left(\frac{i}{q}\right)\right)(-1)^{i}\lambda_{i}\nonumber \\
&=& \lambda_{kq}
\mathscr{L}_n( T_{p|k}\varphi).\nonumber
\end{eqnarray}

Now suppose $k$ is odd. Then manipulating indices and using Eq. (\ref{eq2}) we get 
\begin{eqnarray}
s_{k}&=&
\sum_{i=0}^{q-1}\textstyle
\left(\varphi\left(\frac{kq+q-1-i}{pq}\right)-
\varphi\left(\frac{kq+q-i}{pq}\right)\right)(-1)^{i}\lambda_{kq+q-1-i}\nonumber \\
&=& \lambda_{kq+q-1}\sum_{i=0}^{q-1}\textstyle
\left(T_{p|k}\varphi\left(1-\frac{i+1}{q}\right)-
T_{p|k}\varphi\left(1-\frac{i}{q}\right)\right)(-1)^{i}\lambda_{i}\nonumber\\
&=& \lambda_{kq+q-1}\mathscr{L}_n\left(\overline{T_{p|k}\varphi}\right).\nonumber
\end{eqnarray}

In conclusion,  
\begin{displaymath}\mathscr{L}_{n+1}(\varphi)=\sum_{{k\ \hbox{\tiny  even}}\atop {0\le k<p}}\theta^n(\lambda_k)\mathscr{L}_n(T_{p|k}\varphi)+
\sum_{{k\ \hbox{\tiny  odd}}\atop {0\le k<p}}\theta^n(\lambda_k)\mathscr{L}_n\left(\overline{T_{p|k}\varphi}\right),
\end{displaymath} and the lemma follows. \qed
\end{pf}

\begin{rem}\label{direct_sum}{\em The expression for $S(\mathscr{L}(\varphi))$ given in Lemma \ref{shift_op_on_Lphi} is unique, in the following sense: if $S(\mathscr{L}(\varphi))=\sum_{k=0}^{p-1}\lambda_k \mathbf{v}^{(k)}$   with $\mathbf{v}^{(k)}\in \Lambda_0$, then $\mathbf{v}^{(k)}$ is equal to $\mathscr{L}(T_{p|k}\varphi)$ or $\mathscr{L}\left(\overline{T_{p|k}\varphi}\right)$, according as $k$ is even or odd.  In fact, the $n$th entry  $u_n$ of an element $\mathbf{u}\in \Lambda_0$ is a linear combination of the $\lambda_i$ with $0\le i<p^n$, and equations (\ref{eq1}) and (\ref{eq2}) show that $\theta^n(\lambda_k) u_n$ is a linear combination of the   $\lambda_i$ with $kp^n\le i<(k+1)p^n$. It follows  that  $\sum_{k=0}^{p-1}\lambda_k \Lambda_0$ is a direct sum. 
}\end{rem}

If $\varphi$ is a $p$-fractal, then one can  easily   see that there exists a finite dimensional  subspace of $\Qset^\mathscr{I}$ containing $\varphi$ and stable under the operators $T_{q|b}$ and reflections $\psi\mapsto\overline{\psi}$. For   $\overline{\varphi}$ is also a $p$-fractal, by \cite[Lemma 4.2]{pfractals1}, and therefore  the  subspace $V$ of $\Qset^\mathscr{I}$  spanned by  $\varphi$, $\overline{\varphi}$, and all their transforms under the operators $T_{q|b}$ is finite dimensional,  and evidently stable under the $T_{q|b}$. A simple calculation done in the proof of \cite[Lemma 4.2]{pfractals1} shows that $\overline{T_{q|b}\psi}=T_{q|q-1-b}\overline{\psi}$, for any $\psi$, so $V$ is also stable under reflections.
If  $M=\mathscr{L}(V)$, then $M$ is a finite dimensional subspace  of $\Lambda_0$ containing $\mathscr{L}(\varphi)$, and Lemma \ref{shift_op_on_Lphi} shows that $S(M)\subseteq \sum_{k=0}^{p-1}\lambda_k M$. This is essentially the characterization of $p$-fractals that we were after.

\begin{defn}\label{defn_regular}{\em 
A coherent sequence $\mathbf{u}$ is {\bf regular} if there exists a finite dimensional subspace $M$ of $\Lambda_0$ containing $\mathbf{u}$ with  $S(M)\subseteq \sum_{k=0}^{p-1}\lambda_k M$.
}\end{defn}

\begin{thm}\label{p-fractal_iff_regular} $\varphi:\mathscr{I}\to \Qset$ is a $p$-fractal if and only if $\mathscr{L}(\varphi)$ is regular. 
\end{thm}

\begin{pf} The ``only if'' direction  was shown before Definition   \ref{defn_regular}.  Now  suppose  $\mathscr{L}(\varphi)$ is regular, and let $M\ni \mathscr{L}(\varphi)$ be as in    \ref{defn_regular}.   
By Eq. (\ref{R_and_S}), $S(M+R(M))$ is contained in $\sum_{k=0}^{p-1} \lambda_k (M+R(M))$. So we may replace $M$ by $M+R(M)$; this allows us to assume that $V=\overline{V}$, where $V=\mathscr{L}^{-1}(M)$. Now $V$ contains $\varphi$, and Lemma \ref{onto-coherent} tells us that $V$ is finite dimensional. It remains to show that $V$ is $p$-stable. Suppose $\eta\in V$. Then $\mathscr{L}(\eta)\in M$, and so $S(\mathscr{L}(\eta))=\sum_{k=0}^{p-1}\lambda_k \mathbf{v}^{(k)}$,
with $\mathbf{v}^{(k)}\in M$.  Comparing with    Lemma \ref{shift_op_on_Lphi} and using Remark \ref{direct_sum},  we find that $T_{p|k}\eta$ ($k$ even) and   $\overline{T_{p|k}\eta}$ ($k$ odd) are in $\mathscr{L}^{-1}(M)=V$. Since $V=\overline{V}$, all the $T_{p|k}\eta$  are in $V$, and $V$ is $p$-stable. 
\qed
\end{pf}

\begin{rem}\label{regularity_preserved}{\em
Note that if $S(M)\subseteq \sum_{k=0}^{p-1}\lambda_k M$ and $S(N)\subseteq \sum_{k=0}^{p-1}\lambda_k N$,
then $S(M+N)\subseteq \sum_{k=0}^{p-1}\lambda_k (M+N)$ and $S(MN)\subseteq \sum_{k=0}^{p-1}\lambda_k
MN$. Consequently, if $\mathbf{u}$ and $\mathbf{v}$ are regular, then so are $\mathbf{u}+\mathbf{v}$ and 
$\mathbf{u}\mathbf{v}$.
}\end{rem}

\section{Operations preserving strong rationality}\label{operations}

In   this section we show that  strong rationality is preserved under a number of operations. More precisely, if $f(\mathbf{x})\in \Bbbk\llbracket \mathbf{x}\rrbracket$ and $g(\mathbf{y})\in \Bbbk\llbracket \mathbf{y} \rrbracket$ 
 are strongly rational, then so are powers of these   series, 
 the product $f(\mathbf{x})g(\mathbf{y})$ and the sum 
  $f(\mathbf{x})+g(\mathbf{y})$.

\begin{prop}\label{pfractal-preserved-by-powers}
If $f$ is strongly rational, then so is any power $f^{m}$ of $f$.
\end{prop}

\begin{pf}
Let $V$ be a finite dimensional $p$-stable  subspace of $\Qset^\mathscr{I}$ containing $\varphi_f$ and the constant function 1. Let $\ell(t)=mt$, and define $V^*$ as in Lemma 3.8 of \cite{pfractals1}. Namely, $V^*$ consists of all $\varphi:\mathscr{I}\to \Qset$ whose restriction to each interval $[d/m,(d+1)/m]$, $d=0,\ldots,m-1$, is $t\mapsto \psi(\ell(t)-d)$ for some $\psi\in V$. Then $V^*$ is a finite dimensional $p$-stable  subspace of $\Qset^\mathscr{I}$, by \cite[Lemma 3.8]{pfractals1}, and it contains $\varphi_{f^m}$, since $\varphi_{f^m}(t)=\varphi_f(\ell(t))$ on $[0,1/m]$ and $\varphi_{f^m}(t)=1$ on each $[d/m,(d+1)/m]$ with $d>0$. 
\qed \end{pf}

\begin{prop}\label{pfractal-preserved-by-products}
Suppose  $f\in \Bbbk\llbracket x_1,\ldots,x_s\rrbracket$ and $g \in \Bbbk\llbracket y_1,\ldots,y_{s'}\rrbracket$. 
If $f $ and $g $ are strongly rational, then so is 
$f g $.
\end{prop}

\begin{pf} Let $q$ be an arbitrary power of $p$, and $0\le a\le q$. 
Let $I$ and $J$ be the colon ideals $(\mathbf{x}^q:f^{a})$ and $(\mathbf{y}^q:g^{a})$. 
Then $((\mathbf{x}^q,\mathbf{y}^q):f^{a})=(I,\mathbf{y}^q)$, and 
$((\mathbf{x}^q,\mathbf{y}^q):(fg)^{a})=(I,J)$. Then
\begin{equation}
\begin{array}{rcl}
\deg(\mathbf{x}^q,\mathbf{y}^q,(fg)^{a})&=&q^{s+s'}-\deg(I,J)\\
&=& q^{s+s'}-\deg(I)\deg(J)\\
&=& q^{s+s'}-(q^{s}-\deg(\mathbf{x}^q,f^{a}))(q^{s'}-\deg(\mathbf{y}^q,g^{a}))\\
&=& q^{s'}\deg(\mathbf{x}^q,f^{a})+q^{s}\deg(\mathbf{y}^q,g^{a})-
\deg(\mathbf{x}^q,f^{a})\deg(\mathbf{y}^q,g^{a}),
\end{array}\nonumber
\end{equation}
and dividing by $q^{s+s'}$ we conclude that 
$\varphi_{fg} =\varphi_{f} +
\varphi_{g} -\varphi_{f} \varphi_{g}$.
The result then follows from the fact that the $p$-fractals form a subalgebra of $\Qset^\mathscr{I}$.
\qed \end{pf}

\begin{prop}\label{pfractal-preserved-by-sums}
Suppose $f\in \Bbbk\llbracket x_1,\ldots,x_s\rrbracket$ and $g \in \Bbbk\llbracket y_1,\ldots,y_{s'}\rrbracket$. 
If $f$ and $g$ are strongly rational, then so is 
$f+g$.
\end{prop}

\begin{pf}
In view of Theorem \ref{p-fractal_iff_regular}, it suffices to show that   $\mathscr{L}(\varphi_{f+g})$ is regular whenever $\mathscr{L}(\varphi_f)$ and $\mathscr{L}(\varphi_g)$ are regular. But $\mathscr{L}(\varphi_{f+g})=\mathscr{L}(\varphi_{f})\mathscr{L}(\varphi_{g})$, as noted after Definition \ref{defn-L}, and the result follows from Remark \ref{regularity_preserved}.
\qed
\end{pf}

We proved in \cite{pfractals1} that $\varphi_{f}$ is a $p$-fractal  for any power series  
 $f$ in two variables with coefficients in a finite field. This, combined with the above results,
shows the rationality of the Hilbert-Kunz series of a large family of power series. In particular, we have the following:

\begin{thm} Suppose $f_i$ and $f_{ij}$ are 2 variable power series over a finite field. Then $\sum_i f_i(x_i,y_i)$ and  $\sum_i \prod_j f_{ij}(x_{ij},y_{ij})$ are strongly rational. In particular, the Hilbert-Kunz series of these power series are rational.
\end{thm}

\section{Examples}\label{examples}

In preparation for the examples that will follow, we introduce some notation. $\Delta$ is the element $\mathscr{L}(\varphi_x)=\mathscr{L}(t)$ of $\Lambda$. So $\Delta=(\delta_1,p^{-1}\delta_p,p^{-2}\delta_{p^2},\ldots)$.   Since $\mathscr{L}(1-t)=-\mathscr{L}(t)$, $R(\Delta)=-\Delta$. If $\mathbf{v}=(v_0,v_1,v_2,\ldots) \in \Lambda_0$, then $\mathbf{v}\cdot \Delta=\alpha(v_0)\Delta$, and in particular $\Delta\cdot \Delta=\Delta$. This can be seen by writing $\mathbf{v}$ as $\mathscr{L}(\varphi)$, for some $\varphi$, and  observing that $(-1)^i\lambda_i \delta_q=\delta_q$ for any $i<q$---an easy  consequence of Eq. (\ref{eq5}).

\subsection{Example 1}\label{ex1}

Let $\Bbbk=\Zset/(3)$, $f=y^3-x^4+x^2y^2$ and $g=xy(x+y)$. We shall calculate the Hilbert-Kunz series of $f(x_1,y_1)+f(x_2,y_2)$ and of $f(x_1,y_1)+g(x_2,y_2)$. In what follows we will indicate that two functions $\varphi$ and $\psi$ differ by a linear function by writing $\varphi\approx \psi$. Set $\varphi=\varphi_f$. In the first example of \cite[Section 6]{pfractals1} we found that $9T_{3|0}\varphi\approx \bar{\varphi}$,  $9T_{3|1}\varphi\approx\varphi$ and $9T_{3|2}\varphi$ is a constant function. Set $\mathbf{a}=\mathscr{L}(\varphi)$ and $\mathbf{b}=\mathscr{L}(\bar{\varphi})=R(\mathbf{a})$. Lemma \ref{shift_op_on_Lphi} tells us that $9S(\mathbf{a})=9(\lambda_0\mathscr{L}(T_{3|0}\varphi)+\lambda_1\mathscr{L}(\overline{T_{3|1}\varphi}))=(\mathbf{b}+(\mbox{constant})\cdot \Delta)+\lambda_1(\mathbf{b}+(\mbox{constant})\cdot \Delta)$. Now $9\mathbf{a}$ is the sequence $(9\lambda_0,8\lambda_0-\lambda_1,\ldots)$ while $\mathbf{b}=R(\mathbf{a})=(-\lambda_0,\ldots)$. So $9S(\mathbf{a})-\mathbf{b}-\lambda_1 \mathbf{b}=(9\lambda_0,\ldots)$. It follows that:
\begin{equation}\label{S_of_a}
9S(\mathbf{a})=\mathbf{b}+\lambda_1\mathbf{b}+9\Delta. 
\end{equation}
Applying $R$ and using Eq. (\ref{R_and_S}) of Section \ref{sequences} we find:
\begin{equation}\label{S_of_b}
9S(\mathbf{b})=\lambda_2 \mathbf{a}+\lambda_1\mathbf{a}-9\lambda_2\Delta.
\end{equation}

Set $\psi=\varphi_g$. The results of Example 2 of \cite{pfractals1}  show that $9T_{3|0}\psi\approx\psi$,  $9T_{3|1}\psi\approx\bar{\psi}$ and $9T_{3|2}\psi$ is a constant function. So if we set $\mathbf{c}=\mathscr{L}(\psi)$, and note that $9\mathbf{c}=(9\lambda_0,7\lambda_0-2\lambda_1,\ldots)$, an argument like that of the last paragraph gives:
\begin{equation}\label{S_of_c}
9S(\mathbf{c})=\mathbf{c}+\lambda_1\mathbf{c}+(6\lambda_0-3\lambda_1)\Delta. 
\end{equation}

To proceed further we define a $\Qset\llbracket {\mathsf{z}}\rrbracket$-valued bilinear function $r$. For sequences $\mathbf{u}=(u_0,u_1,\ldots)$ and $\mathbf{v}=(v_0,v_1,\ldots)$  in $\Lambda$, set
\begin{equation}
r(\mathbf{u},\mathbf{v})=(1-27{\mathsf{z}})\cdot \sum_{n=0}^{\infty} \alpha(u_nv_n)(81{\mathsf{z}})^n.\nonumber
\end{equation}
The remarks made after Definition \ref{defn-L} tell us that  the Hilbert-Kunz series of $f(x_1,y_1)+f(x_2,y_2)$ is $(1-27{\mathsf{z}})^{-1}\cdot r(\mathbf{a},\mathbf{a})$. Similarly, the Hilbert-Kunz series of $f(x_1,y_1)+g(x_2,y_2)$ is $(1-27{\mathsf{z}})^{-1}\cdot r(\mathbf{a},\mathbf{c})$. The bilinear function $r$ has the following basic properties:
\begin{itemize}
\item[(A)] $r(\mathbf{u},\mathbf{v})=(1-27{\mathsf{z}})\cdot \alpha(u_0v_0)+{\mathsf{z}}\cdot r(9S(\mathbf{u}),9S(\mathbf{v}))$;
\item[(B)] $r(R(\mathbf{u}),R(\mathbf{v}))=r(\mathbf{u},\mathbf{v})$;
\item[(C)] If  $\mathbf{u}$ and $\mathbf{v}$ satisfy condition (1) of Definition \ref{defn-coherent}, then $r(\lambda_i \mathbf{u},\lambda_j \mathbf{v})=r(\mathbf{u},\mathbf{v})$ if $i=j$, and  $r(\lambda_i \mathbf{u},\lambda_j \mathbf{v})=0$ otherwise;
\item[(D)] If $\mathbf{u}\in \Lambda_0$, $r(\mathbf{u},\Delta)$ is the constant $\alpha(u_0)$. (To see this, note that since $\mathbf{u}\cdot \Delta=\alpha(u_0)\Delta$, $\alpha(u_n\Delta_n)=\alpha(u_0)\alpha(\Delta_n)=3^{-n}\alpha(u_0)$.)
\end{itemize}

Using (A)--(D), Eq. (\ref{S_of_a}) and the fact that $\alpha(a_0)=1$ we find:
\begin{eqnarray}
r(\mathbf{a},\mathbf{a})&=& 1-27{\mathsf{z}}+{\mathsf{z}}\cdot r(9S(\mathbf{a}),9S(\mathbf{a})) \nonumber\\
&=&  1-27{\mathsf{z}}+{\mathsf{z}}\cdot r(\mathbf{b}+\lambda_1\mathbf{b}+9\Delta,\mathbf{b}+\lambda_1\mathbf{b}+9\Delta) \nonumber\\
&=&  1-27{\mathsf{z}}+{\mathsf{z}}(2r(\mathbf{b},\mathbf{b})-9-9+81) \nonumber\\
&=& 1+36{\mathsf{z}}+2{\mathsf{z}}\cdot r(\mathbf{a},\mathbf{a}).\nonumber
\end{eqnarray}
So $r(\mathbf{a},\mathbf{a})=(1+36 {\mathsf{z}})/(1-2{\mathsf{z}})$, and the Hilbert-Kunz series of $h=f(x_1,y_1)+f(x_2,y_2)$ is $(1+36 {\mathsf{z}})/((1-2{\mathsf{z}})(1-27{\mathsf{z}}))$. From this rational description of $\hks(h)$ we can obtain the Hilbert-Kunz multiplicity $\mu$ of $h$  as follows. Since $e_n(h)=\mu \cdot 27^n+O(9^n)$,  $\sum_{n=0}^\infty (e_n(h)-\mu\cdot 27^n){\mathsf{z}}^n$ converges on a
neighborhood of $1/27$. So  $\hks(h)-\mu/(1-27{\mathsf{z}})$ is holomorphic on this neighborhood, and 
 $\mu= \lim_{{\mathsf{z}} \to 1/27} (1-27{\mathsf{z}})\hks(h)=(1+36/27)/(1-2/27)=63/25$.

Turning to $f(x_1,y_1)+g(x_2,y_2)$, we compute $r(\mathbf{a},\mathbf{c})$ and $r(\mathbf{b},\mathbf{c})$ in similar fashion. Noting that $\alpha(b_0)=-1$ and $\alpha(c_0)=1$ we find:
\begin{eqnarray}
r(\mathbf{a},\mathbf{c})&=& 1-27{\mathsf{z}}+{\mathsf{z}}\cdot r(9S(\mathbf{a}),9S(\mathbf{c})) \nonumber\\
&=&  1-27{\mathsf{z}}+{\mathsf{z}}\cdot r(\mathbf{b}+\lambda_1\mathbf{b}+9\Delta,\mathbf{c}+\lambda_1\mathbf{c}+6\Delta-3\lambda_1 \Delta) \nonumber\\
&=&  1-27{\mathsf{z}}+{\mathsf{z}}(2r(\mathbf{b},\mathbf{c})-6+3+9+54) \nonumber\\
&=& 1+33{\mathsf{z}}+2{\mathsf{z}}\cdot r(\mathbf{b},\mathbf{c});\nonumber
\end{eqnarray}
\begin{eqnarray}
r(\mathbf{b},\mathbf{c})&=& -1+27{\mathsf{z}}+{\mathsf{z}}\cdot r(9S(\mathbf{b}),9S(\mathbf{c})) \nonumber\\
&=&  -1+27{\mathsf{z}}+{\mathsf{z}}\cdot r(\lambda_2\mathbf{a}+\lambda_1\mathbf{a}-9\lambda_2\Delta,\mathbf{c}+\lambda_1\mathbf{c}+6\Delta-3\lambda_1 \Delta) \nonumber\\
&=&  -1+27{\mathsf{z}}+{\mathsf{z}}(r(\mathbf{a},\mathbf{c})-3) \nonumber\\
&=& -1+24{\mathsf{z}}+{\mathsf{z}}\cdot r(\mathbf{a},\mathbf{c}).\nonumber
\end{eqnarray}
The above equations  yield  $r(\mathbf{a},\mathbf{c})=(1+31{\mathsf{z}}+48 {\mathsf{z}}^2)/(1-2{\mathsf{z}}^2)$, and the Hilbert-Kunz series of $f(x_1,y_1)+g(x_2,y_2)$ is $(1+31{\mathsf{z}}+48 {\mathsf{z}}^2)/((1-2{\mathsf{z}}^2)(1-27{\mathsf{z}}))$. In particular, the Hilbert-Kunz multiplicity   is $(1+31/27+48/729)/(1-2/729)=1614/727$.

\bigskip
\bigskip

Our next example uses  further  notation and results from Example 3 of  \cite{pfractals1}. Let $\lambda\in \Bbbk\setminus\{0,1\}$, and $h=\prod_{i=1}^4 h_i$, where the $h_i$ are $x$, $y$, $x+y$ and $x+\lambda y$.  
Set $B=\Bbbk\llbracket x,y \rrbracket/(h)$. 
   $X_2=X_2(h)$ is the set of equivalence  classes of  two-generator  finite colength   ideals of $B$. 
   In  $X_2$ we have:
     
\begin{enumerate}
\item $O$, $E_1$, $E_\infty$ and $E_0$; the   classes containing the ideals $(1)$, $(h_1h_2,h_3h_4)$, 
$(h_1h_3,h_2h_4)$ and $(h_1h_4,h_2h_3)$. (In the language of \cite{pfractals1}, these are ``reflections'' of $O$.) 

\item $E_\beta$, with $\beta\in \Bbbk\setminus\{0,1\}$; the class ``represented by $P_\beta$'', in the language of Proposition 6.3 of \cite{pfractals1}. Explicitly, $E_\lambda$ is the class of $(x,y)$, while if $\beta\ne \lambda$, $E_\beta$ is the class of $((x+y)^2,(x+\lambda y)(\lambda x+\beta y))$.

\end{enumerate}

In \cite{pfractals1} we defined an involution $R:X_2\to X_2$, and noted that $R$ fixes all the $E_\beta$, $\beta\in \Bbbk\cup \{\infty\}$, except for $E_\lambda$. Now let $Y$ consist of all the $E_\beta$, together with $O$ and $R(E_\lambda)$, so that $Y$ is stable under $R$. In \cite{pfractals1} we also introduced ``magnification operators'' $\tau_i=\tau_{p|(i,i,i,i)}:X_2\to X_2$, $0\le i<p$.  It is  not hard to see that $Y$ is stable under these operators---see \cite{paul3}. 

Now to each $C$ in $Y$ we attach a $\varphi_C:\mathscr{I}\to \Qset$ (and an $\mathbf{e}_C=\mathscr{L}(\varphi_C)$ in $\Lambda_0$) as follows. In \cite{pfractals1} we attached a function $\varphi_I:\mathscr{I}\to \Qset$ to each ideal $I$ in $\Bbbk\llbracket x,y\rrbracket$ whose image in $B$ has finite colength; namely $\varphi_I(a/q)=q^{-2}\deg(I^{[q]},h^a)$. When  $C=E_\lambda$ (resp. $R(E_\lambda)$) we  set  $\varphi_C=\varphi_{(x,y)}$  (resp. $\overline{\varphi_{(x,y)}}$).  For any other $C$ we choose an  ideal $I\supset ( h)$ of $\Bbbk\llbracket x,y\rrbracket$ whose image in $B$ lies in $C$, and we set $\varphi_C(t)=\varphi_I(t)-(\deg I)t$. This $\varphi_C$ is independent of the choice of  $I$, and $\varphi_C(0)=\varphi_C(1)=0$. Using the fact that $R(C)=C$ we find that $\overline{\varphi_C}=\varphi_C$, so that $R(\mathbf{e}_C)=\mathbf{e}_C$. 

Suppose now  that $C$ is some $E_\beta$.  Let $\varphi_i$ be the function attached as above to $\tau_i(C)\in Y$. Set $\mathbf{e}_i=\mathscr{L}(\varphi_i)$ for even $i$, and $\mathbf{e}_i=\mathscr{L}(\overline{\varphi_i})$ for odd $i$.  (Note that when $\beta=\lambda$, $\varphi_C=\varphi_{(x,y)}$ is constant on $[1/2,1]$, so that $\varphi_i=0$ for $i\ge p/2$.) Lemma 3.6 of \cite{pfractals1} shows that $T_{p|i}(p^2\varphi_C)\approx \varphi_i$. Invoking Lemma \ref{shift_op_on_Lphi} we find that 
\begin{equation}\label{eq5.1}
p^2S(\mathbf{e}_C)=\sum \lambda_i\mathbf{e}_i+\sum (-1)^i a_i\lambda_i\Delta,
\end{equation}
 where the $a_i\in \Qset$ and the sums run over all $i<p$ when $\beta\ne \lambda$, and all $i<p/2$ when $\beta=\lambda$. 

To completely determine the shift rules, we need to calculate the $a_i$. This is straightforward when $\beta=0$, $1$ or $\infty$. (In these cases, $\varphi_C(t)=4(t-t^2)$ and $\tau_i(C)=C$, for all $i$.) For other $\beta$, it can be shown that  the $a_i$ can be expressed in terms of certain ``syzygy gaps''. In \cite{paul4}, the first author shows that each of these syzygy gaps is 0 or 2, a result that yields the following theorem, whose proof we omit:

\begin{thm}\label{thm5.1}
If $p=2m+1$, then $a_i=8(m-i)$. 
\end{thm}

\subsection{Example 2}

Suppose $\Bbbk=\Zset/(7)$, $f=x^3+y^4$ and $g=x^4+xy^3$. We shall calculate the Hilbert-Kunz series of $f(x_1,y_1)+g(x_2,y_2)$. Let $\zeta$, $\eta$, $\varphi$ and $\psi$ be $\varphi_{x^3}$, $\varphi_{y^4}$, $\varphi_{f}$ and $\varphi_g$, respectively. Set $\mathbf{u}=\mathscr{L}(\varphi)$, so that $\mathbf{u}=\mathscr{L}(\zeta)\mathscr{L}(\eta)$. 
We begin by working out the shift rule for $\mathbf{u}$. 

Evidently $\zeta(t)=\min(3t,1)$. So $7T_{7|2} \zeta\approx \zeta$, while $T_{7|i}\zeta$ is
 linear for all other $i$. Lemma \ref{shift_op_on_Lphi}  then tells us that $7S(\mathscr{L}(\zeta))=\lambda_2\mathscr{L}(\zeta)+ $a linear combination of the $\lambda_i\Delta$. Now $7\mathscr{L}(\zeta)$ is the sequence $(7\lambda_0,3\lambda_0-3\lambda_1+\lambda_2,\ldots)$. We conclude:
\begin{equation}
7S(\mathscr{L}(\zeta))=\lambda_2\mathscr{L}(\zeta)+(3\lambda_0-3\lambda_1)\Delta.\nonumber
\end{equation}

Similarly $\eta(t)=\min(4t,1)$, and  so $7T_{7|1} \eta\approx\bar{\eta}$, while $T_{7|i}\eta$ is
 linear for all other $i$. Then $7S(\mathscr{L}(\eta))=\lambda_1\mathscr{L}(\eta)+$a linear combination of the $\lambda_i\Delta$. Since $7\mathscr{L}(\eta)=(7\lambda_0,4\lambda_0-3\lambda_1,\ldots)$, we find:
\begin{equation}
7S(\mathscr{L}(\eta))=\lambda_1\mathscr{L}(\eta)+(4\lambda_0-4\lambda_1)\Delta.\nonumber
\end{equation}
  Multiplying the last two displayed formulas, using Eq. (\ref{eq4}) and noting that $\Delta\cdot\Delta=\mathscr{L}(\zeta)\cdot \Delta=\mathscr{L}(\eta)\cdot \Delta=\Delta$, we find:
  \begin{equation}\label{shiftrule1}
  49S(\mathbf{u})=(\lambda_1+\lambda_2+\lambda_3)\mathbf{u}+(21\lambda_0-16\lambda_1+9\lambda_2-4\lambda_3)\Delta.
  \end{equation}

Now set $\mathbf{a}=\mathscr{L}(\psi)$ and $\mathbf{b}=R(\mathbf{a})$. We shall work out the shift rules for $\mathbf{a}$ and $\mathbf{b}$.  A linear change of variables takes $g$ into $h=xy(x+y)(x+3y)$, and we are in the situation described in the paragraphs preceding this example.  Let $M\in X_2(h)$ be the ideal class $E_3$ of the ideal $(x,y)$. Let $M^*=R(E_3)$ and $D=\tau_2(E_3)$. Then $\mathbf{a}=\mathbf{e}_M$ and $\mathbf{b}=\mathbf{e}_{M^*}$; let $\mathbf{d}=\mathbf{e}_D$.  

One can show that the $\tau_i(M)$, $0\le i\le 6$, are the classes $M$, $M^*$, $D$, $M$, $O$, $O$ and $O$, respectively. Eq. (\ref{eq5.1}) and Theorem \ref{thm5.1} then give us:
\begin{equation}\label{shiftrule2}
49 S(\mathbf{a})=(\lambda_0+\lambda_1)\mathbf{a}+\lambda_3\mathbf{b}+\lambda_2\mathbf{d}+(3\lambda_0-2\lambda_1+\lambda_2)8\Delta.
\end{equation}
Applying $R$ we find:
\begin{equation}\label{shiftrule3}
49 S(\mathbf{b})=\lambda_3\mathbf{a}+(\lambda_5+\lambda_6)\mathbf{b}+\lambda_4\mathbf{d}+(-\lambda_4+2\lambda_5-3\lambda_6)8\Delta.
\end{equation}
Furthermore it may be shown that the $\tau_i(D)$, $0\le i\le 6$, are the classes $D$, $M$, $M^*$, $D$, $M$, $M^*$ and  $D$, respectively. Eq. (\ref{eq5.1}) and Theorem \ref{thm5.1}  then give the shift rule:
\begin{eqnarray}\label{shiftrule4}
49 S(\mathbf{d})&=&(\lambda_4+\lambda_5)\mathbf{a}+(\lambda_1+\lambda_2)\mathbf{b}+(\lambda_0+\lambda_3+\lambda_6)\mathbf{d}\nonumber\\
&&+(3\lambda_0-2\lambda_1+\lambda_2-\lambda_4+2\lambda_5-3\lambda_6)8\Delta.
\end{eqnarray}

If $\mathbf{v}$ and $\mathbf{v}'$ are in $\Lambda$, set 
\begin{equation}
r(\mathbf{v},\mathbf{v}')=(1-343{\mathsf{z}})\cdot \sum_{n=0}^\infty\alpha(v_nv'_n)(2401{\mathsf{z}})^n.\nonumber
\end{equation}
 Results analogous to (A)--(D) of Example 1 evidently hold, and the Hilbert-Kunz series of $f(x_1,y_1)+g(x_2,y_2)$ is $(1-343{\mathsf{z}})^{-1}\cdot r(\mathbf{u},\mathbf{a})$, since $\mathbf{u}=\mathscr{L}(\varphi_{f})$ and $\mathbf{a}=\mathscr{L}(\varphi_g)$.  Note that $\alpha(u_0)$, $\alpha(a_0)$, $\alpha(b_0)$ and $\alpha(d_0)$ are 1, 1, $-1$ and 0, respectively. ($R(\mathbf{d})=\mathbf{d}$ gives the last of these facts.)  As in Example 1 we find:
\begin{eqnarray}
r(\mathbf{u},\mathbf{a}) & = & 1-343{\mathsf{z}}+{\mathsf{z}}\cdot r(49S(\mathbf{u}),49 S(\mathbf{a}))\nonumber\\
r(\mathbf{u},\mathbf{b}) & = & -1+343{\mathsf{z}}+{\mathsf{z}}\cdot r(49S(\mathbf{u}),49 S(\mathbf{b}))\nonumber\\
r(\mathbf{u},\mathbf{d}) & = & {\mathsf{z}}\cdot r(49S(\mathbf{u}),49 S(\mathbf{d}))\nonumber
\end{eqnarray}
Using the shift rules (\ref{shiftrule1})--(\ref{shiftrule4}) and arguing as in Example 1 we get:
\begin{eqnarray}
r(\mathbf{u},\mathbf{a}) & = & 1+490{\mathsf{z}}+{\mathsf{z}}(r(\mathbf{u},\mathbf{a})+r(\mathbf{u},\mathbf{b})+r(\mathbf{u},\mathbf{d}))\nonumber\\
r(\mathbf{u},\mathbf{b}) & = & -1+339{\mathsf{z}}+{\mathsf{z}}\cdot r(\mathbf{u},\mathbf{a})\nonumber\\
r(\mathbf{u},\mathbf{d}) & = &831{\mathsf{z}}+{\mathsf{z}}(2r(\mathbf{u},\mathbf{b})+r(\mathbf{u},\mathbf{d}))\nonumber
\end{eqnarray}
The above system gives $r(\mathbf{u},\mathbf{a})=(1+488{\mathsf{z}}+679{\mathsf{z}}^2+339{\mathsf{z}}^3)/(1-2{\mathsf{z}}-{\mathsf{z}}^3)$. We conclude that the Hilbert-Kunz series of $f(x_1,y_1)+g(x_2,y_2)$ is  \begin{equation}
\frac{1+488{\mathsf{z}}+679{\mathsf{z}}^2+339{\mathsf{z}}^3}{(1-343{\mathsf{z}})(1-2{\mathsf{z}}-{\mathsf{z}}^3)}.\nonumber
\end{equation}

\section{The Hilbert-Kunz function of $z^D-h(x,y)$}\label{application}
 
 We next investigate in more detail the Hilbert-Kunz function of $g=z^D-h(x,y)$, where $h\ne 0$ is in the maximal ideal of $\Bbbk\llbracket x,y\rrbracket$. We allow $\Bbbk$ to be infinite. Write $D=p^cE$, with $\gcd(E,p)=1$. A  key result of this section, Theorem \ref{6.6}, is an estimate in a (one sided) neighborhood of $1/E$ for a certain function $\psi:\mathscr{I}\to \Qset$ constructed from $h$.   Then we will see how to express the  $e_n(g)$ in terms of values of $\psi$ at two sequences of  points, each of which  approaches $1/E$ as $n\to \infty$. We shall make use of a contraction operator having $1/E$ as its fixed point, together with ideal class ideas from \cite{pfractals1}. 
 
 We may write $h=\prod_{i=1}^rh_i^{d_i}$, with the $h_i$ pairwise prime irreducibles. Let $G=\prod_{i=1}^r h_i$, and $X_2(G)$ be as in \cite{pfractals1}. Choose a constant $K$ such that every ideal class $C$ in $X_2(G)$ contains an ideal which is the image of some $I\subseteq \Bbbk\llbracket x,y\rrbracket$ with $\deg I\le K$. Remark 3.3 of \cite{pfractals1} shows this is possible. Given $C$ we fix such an $I=I_C$ and let $\Phi_C$ be the function $\mathscr{I}^r \to \Qset$ attached to $I$, i.e., the function $(a_1/q,\ldots,a_r/q)\mapsto q^{-2}\deg(I^{[q]},\prod_{i=1}^r h_i^{a_i})$.  Note that $|\Phi_C|\le \deg I\le K$. If $C$ is the trivial class, $\Phi_C$ is linear, while if $C$ is the class of $(x,y)$, then $\Phi_C\approx\Phi_G$, where $\Phi_G(a_1/q,\ldots,a_r/q)=q^{-2}\deg(x^q,y^q,\prod_{i=1}^r h_i^{a_i})$.
 
 \begin{defn}\label{6.1}{\em
 Suppose $a=(a_1,\ldots,a_r)$ with $0\le a_i<d_i$. Let $X(a)$ consist of all $t=(t_1,\ldots, t_r)$ in $\mathscr{I}^r$ with each $d_it_i$ in $[a_i,a_{i}+1]$; the $X(a)$ cover $\mathscr{I}^r$. $\Phi:\mathscr{I}^r\to \Qset$ is of {\bf $h$-type} if for each $X(a)$ there is a $C$ in $X_2(G)$ such that the restriction of $\Phi$ to $X(a)$ is $t\mapsto \Phi_C(dt-a)+(\text{linear})$, where $dt$ denotes the component-wise product of the vectors $d=(d_1,\ldots,d_r)$ and $t$. 
 }\end{defn}
 
 Note that $\Phi_h:\mathscr{I}^r\to \Qset$ with   $\Phi_h(c_1/q,\ldots,c_r/q)=q^{-2}\deg(x^q,y^q,\prod_{i=1}^r h_i^{d_ic_i})$ is of $h$-type. For if $a=(0,\ldots,0)$, then on $X(a)$,    $\Phi_h(t)=\Phi_G(dt)$, while $\Phi_h$ is the constant function  1 on all other $X(a)$. 
 
 \begin{lem}\label{6.2}If $\Phi$ is of $h$-type, then for all $b=(b_1,\ldots,b_r)$ with $0\le b_i<p$, $\Psi=p^2T_{p|b}\Phi$ is of $h$-type. 
 \end{lem}
 
 \begin{pf}
 What follows is essentially contained in Lemma 3.8 of \cite{pfractals1}. Fix $a$ as in Definition \ref{6.1} and write $a+db$ as $pa^*+b^*$, with $0\le b_i^*<p$. Note that if $t$ is in $X(a)$, then
 $(t+b)/p$ is in $X(a^*)$. 
                   Choose $C$ so that $\Phi(u)\approx\Phi_C(du-a^*) $ on $X(a^*)$.  So 
 $\Phi((t+b)/p)\approx\Phi_C(d(t+b)/p-a^*)$ on $X(a)$.  Multiplying by $p^2$ we find that on $X(a)$, $\Psi(t)\approx p^2\Phi_C((dt-a+b^*)/p)= p^2T_{p|b^*}\Phi_C(dt-a)$. But \cite{pfractals1} tells us that $p^2T_{p|b^*}\Phi_C\approx\Phi_{C'}$, where $C'=\tau_{p|b^*}C$.  Since this result holds for all $a$, $\Psi$ is of $h$-type. 
 \qed
 \end{pf}
 
 Returning to the first paragraph of this section, we write $D=p^cE$ with $\gcd(E,p)=1$. We shall assume $E\ne 1$. Let $\alpha$ be the order of $p$ in $(\Zset/(E))^*$. Set $P=p^\alpha$ and write $P-1$ as $E\gamma$. 
 
 \begin{defn}\label{6.3}{\em 
 $\Psi=p^{2c}T_{p^c|(0,\ldots,0)}\Phi_h$ and  $\Psi_n=\mathscr{J}^n(\Psi)$, where $\mathscr{J}$ is the operator $P^2T_{P|(\gamma,\ldots,\gamma)}$. 
 }\end{defn}
 
 Let $\psi$ and $\psi_n:\mathscr{I}\to \Qset$ be the compositions  of $\Psi$ and $\Psi_n$ with  the diagonal map $t\mapsto (t,\ldots,t)$. If we take $\varphi_h:\mathscr{I}\to \Qset$ to be the function $a/q\mapsto q^{-2}\deg(x^q,y^q,h^a)$, then $\varphi_h$ is the composition of $\Phi_h$ and the diagonal map. It follows that $\psi=p^{2c}T_{p^c|0}\varphi_h$, and that $\psi_n=J^n(\psi)$, where $J$ is the operator $P^2T_{P|\gamma}$. 
 
 Note that the map $[0,1]\to [0,1]$, $w\mapsto (w+\gamma)/P$, stabilizes $\mathscr{I}$, has slope $1/P$ and fixes $E^{-1}$. So this map takes $E^{-1}+w$ to $E^{-1}+w/P$. Our plan is to use this contraction mapping to study the function $\psi$ in a right (resp. left) neighborhood of $E^{-1}$ in $\mathscr{I}$. 
 
 Lemma \ref{6.2} shows  that each $\Psi_n$ is of $h$-type. So the restriction of $\Psi_n$ to any $X(a)$ has distance $\le K$ from a linear function in the uniform metric (since each  $|\Phi_C|\le K$). Choose an open right neighborhood $(0,\epsilon)$ of 0 such that all points $(E^{-1}+w,\ldots,E^{-1}+w)$, with $w$ in this interval and $E^{-1}+w$ in $\mathscr{I}$, lie in a single $X(a)$. Then there exist $\lambda_n$ and $u_n$ such that for all such $w$, $|\Psi_n(E^{-1}+w,\ldots,E^{-1}+w)-\lambda_nw-u_n|\le K$. In other words:
 
 \begin{lem}\label{6.4}
 $|\psi_n(E^{-1}+w)-\lambda_nw-u_n|\le K$ whenever $w$ is small and positive and $E^{-1}+w$ is in $\mathscr{I}$. 
 \end{lem} 
 
 \begin{rem}\label{6.4.1}{\em 
 One gets a similar result for small negative $w$, but it may not be possible to choose the same constants $\lambda_n$ and $u_n$. When no $d_i$ is a multiple of $E$, however, the point $(E^{-1},\ldots,E^{-1})$ is an interior point of an $X(a)$, and so Lemma \ref{6.4} holds for all $w$ in a small 2-sided neighborhood of 0. 
 }\end{rem}
 
 \begin{lem}\label{6.5}
 Suppose we are in the situation of Lemma \ref{6.4}. Then there are real numbers $\lambda$ and $u$ such that $\lambda_n=\lambda P^n+O(1)$ and  $u_n=uP^{2n}+O(1)$. 
 \end{lem}
  
  \begin{pf}
  Replacing $w$ by $w/P$ in Lemma \ref{6.4} and multiplying by $P^2$ (noting that $\psi_{n+1}(E^{-1}+w)=P^2\psi_n(E^{-1}+w/P)$, since $\psi_{n+1}=J(\psi_n)$), we find:
  \begin{equation}
|\psi_{n+1}(E^{-1}+w)-P\lambda_nw-P^2u_n|\le P^2 K.
  \nonumber
  \end{equation}
  Combining this inequality with the inequality obtained by replacing $n$ by $n+1$ in Lemma \ref{6.4} we find that $|(\lambda_{n+1}-P\lambda_n)w+(u_{n+1}-P^2u_n)|\le P^2K+K$. Since this holds for all $w$ in a fixed interval, $\lambda_{n+1}-P\lambda_n$ is $O(1)$, and consequently $u_{n+1}-P^2u_n$ is $O(1)$. Then $\lambda_{n+1}/P^{n+1}-\lambda_n/P^n$ is $O(P^{-n})$; we conclude that $\lambda_n/P^n$ converges to some $\lambda$, and that $\lambda_n/P^n=\lambda+O(P^{-n})$. This gives the result for $\lambda_n$; that for $u_n$ is proved similarly. 
  \qed
  \end{pf}
  
  We now derive the desired estimate for $\psi$ in a neighborhood of $E^{-1}$. 
  
  \begin{thm}\label{6.6}
  Let $u$ and $\lambda$ be as in Lemma \ref{6.5}. Then for $w>0$, with $E^{-1}+w$ in $\mathscr{I}$, we have $\psi(E^{-1}+w)=u+\lambda w+O(w^2)$. 
  \end{thm}
  
\begin{pf}
Fix $\epsilon>0$ so that the result of Lemma \ref{6.4} holds for $w$ in $(0,\epsilon)$. If $w$ is in $(0,\epsilon)$, choose $n$ so that $P^nw$ is in $[\epsilon/P,\epsilon)$. Then $\psi(E^{-1}+w)=P^{-2n}\psi_n(E^{-1}+P^nw)$. Now Lemma \ref{6.4}, with $w$ replaced by $P^nw$, together with Lemma \ref{6.5}, tells us that $\psi_n(E^{-1}+P^nw)-P^{2n}\lambda w-P^{2n}u$ is bounded independently of $w$. Thus $|\psi(E^{-1}+w)-\lambda w-u|\le (\text{constant})P^{-2n}$. Since $P^{-2n}=w^2(P^nw)^{-2}$, and $|P^nw|\ge \epsilon/P$, we get the result. 
\qed\end{pf}  

\begin{rem}\label{6.6.1}{\em
There is a similar result  when $w<0$, with a possibly different $\lambda$. When $E$ does not divide any $d_i$, Remark \ref{6.4.1} shows that Theorem \ref{6.6} holds for $w$ in  a 2-sided neighborhood of 0.
}\end{rem}
  
We now estimate $e_n(g)$, where $g=z^D-h$.  Recall that $D=p^cE$, with $\gcd(E,p)=1$. If $n\ge c$ write $p^{n-c}=s_nE+r_n$, with $0\le r_n<E$. Then  $q=p^n=s_nD+p^cr_n$.  Note that $p^cs_n/q$ is close to $E^{-1}$ for large $n$. Observe in fact:
\begin{equation}
\frac{p^cs_n}{q}=\frac{p^c}{q}\left(\frac{q-p^cr_n}{D}\right)=\frac{1}{E}-\frac{p^{2c}r_n}{Dq}
\nonumber
\end{equation}  
\begin{equation}
\frac{p^c(s_n+1)}{q}=\frac{1}{E}-\frac{p^{2c}r_n}{Dq}+\frac{p^cD}{Dq}=\frac{1}{E}+\frac{p^{2c}(E-r_n)}{Dq}
\nonumber
\end{equation}  

\begin{lem}\label{7.1}
Suppose that $E>1$. Then:
\begin{enumerate}
\item $\deg(x^q,y^q,h^{s_n})=(up^{-2c})q^2-(\lambda^-r_n/D)q+O(1)$, for some $u$ and $\lambda^-$. 
\item $\deg(x^q,y^q,h^{s_n+1})=(up^{-2c})q^2+(\lambda^+r_n/D)q+O(1)$, for some  $\lambda^+$. 
\end{enumerate}
If $E$ does not divide any $d_i$, then $\lambda^+=\lambda^-$. 
\end{lem}

\begin{pf} $\deg(x^q,y^q,h^{s_n})=q^2\varphi_h(s_n/q)=p^{-2c}q^2\psi(p^cs_n/q)$. The observation above, combined with Remark \ref{6.6.1}, shows that $\psi(p^cs_n/q)=u-\lambda^-p^{2c}r_n/(Dq)+O(q^{-2})$, for some $u$ and $\lambda^-$. Multiplying by $p^{-2c}q^2$ we get (1). A similar argument, using Theorem \ref{6.6} itself, shows that $\deg(x^q,y^q,h^{s_n+1})=(u^*p^{-2c})q^2+(\lambda^+r_n/D)q+O(1)$, for some $u^*$ and $\lambda^+$. It is easy to see that  $\deg(x^q,y^q,h^{s_n+1})-\deg(x^q,y^q,h^{s_n})$ is $O(q)$, and it follows that $u^*=u$. When $E$ does not divide any $d_i$ we use Remark \ref{6.6.1} to see that $\lambda^+=\lambda^-$. 
\qed\end{pf}

We now estimate $e_n=e_n(g)$ for $n\ge c$. An easy calculation shows that $e_n=p^c((E-r_n)\deg(x^q,y^q,h^{s_n})+r_n\deg(x^q,y^q,h^{s_n+1}))$.  When $E=1$, the $r_n$ are all 0. So $e_n=D\deg(x^q,y^q,h^{q/D})$, and $e_{n+1}=p^2e_n$ for $n\ge c$. Suppose $E>1$. Lemma \ref{7.1} then gives:
\begin{equation}
e_n=(p^{-2c}Du)q^2-\frac{(\lambda^--\lambda^+)r_n(E-r_n)}{E}\cdot q +O(1).\nonumber
\end{equation}
We have proved:
\begin{thm}\label{7.2}
Let $r_n$ be the remainder when $p^{n-c}$ is divided by $E$. Then there exist $\mu$ and $\mu_1$ such that $e_n(g)=\mu p^{2n}-\mu_1r_n(E-r_n)p^n+O(1)$. If $E=1$ or $E$ does not divide any $d_i$ then $\mu_1=0$. 
\end{thm}

When $\Bbbk$ is finite one can say more:

\begin{thm}\label{7.3}
Suppose $\Bbbk$ is finite. Then $\mu$ and $\mu_1$ are rational and the $O(1)$ term in Theorem \ref{7.2} is eventually periodic. 
\end{thm}
\begin{pf} This is immediate from Theorem \ref{7.2} and the rationality of the Hilbert-Kunz series of $g$, proved earlier. One can also give  the argument sketched next, which avoids the heavy machinery used in proving rationality. Since $\Bbbk$ is finite, $X_2(G)$ is a finite set, according to \cite{pfractals1}. The functions $\Psi_n$ of Definition \ref{6.3} are all of $h$-type. Fix $a=(a_1,\ldots,a_r)$. Then there must exist $M$ and $L$, $L>0$, with $\Psi_M\approx \Psi_{M+L}$ on $X(a)$. Composing with the diagonal map one finds that $\psi_M\approx \psi_{M+L}$ on a right neighborhood of $1/E$ and  on a  left  neighborhood of $1/E$. The information obtained in this way about $\psi$, together with the formulas for $\deg(x^q,y^q,h^{s_n})$ and $\deg(x^q,y^q,h^{s_n+1})$ in terms of values of $\psi$, leads quickly to the desired results.
\qed\end{pf}

The $p^n$ term in Theorem \ref{7.2} can occur. Here are two examples. (The first, well-known, goes back to Kunz.)

\begin{enumerate}
\item Suppose $p\ne 5$ and $g=z^5-x^5y^4$, so that $E=5$ and the $d_i$ are 5 and 4. Then $\mu=5$ and $\mu_1=1/5$. 
\item Suppose $p=7$ and $g=z^{14}-x^6y^6(x^2-y^2)$, so that $E=2$ and the $d_i$ are 6, 6, 1, and 1. It can be shown that for $n\ge 2$, $e_n=(74/7)49^n-6\cdot 7^n-42$. So $\mu=74/7$ and $\mu_1=6$. 
\end{enumerate}

We conclude with a few remarks. Brenner \cite{brenner} has results similar in character to Theorems \ref{7.2} and \ref{7.3}. He requires homogeneity and normality, and so he never gets a $p^n$ term. On the other hand, he can treat the homogeneous coordinate ring of an arbitrary non-singular curve. (Also, the ideal used to define his Hilbert-Kunz functions need not be the homogeneous maximal ideal.) In order to get results like Theorem \ref{7.3}, Brenner, like us, has to assume that  $\Bbbk$ is finite. Whether this condition is really needed is unknown. 

Our results, Brenner's results, and the results of Huneke-McDermott-Monsky \cite{HMM}, suggest that analogues of Theorems \ref{7.2} and \ref{7.3} may hold for  general Hilbert-Kunz functions in dimension 2. \cite{HMM} gives a result like Theorem \ref{7.2} in a general setting, but assuming normality, so that no $p^n$ term occurs. Very little is known when the normality assumption is dropped.  The situation is of course completely different in higher dimensions.

\end{document}